\documentclass[12pt]{amsart}
\usepackage{amssymb,amsmath}
\usepackage{MnSymbol}
 \usepackage{tikz}
\numberwithin{equation}{section}

\def\di{\partial}
\def\dib{\bar\partial}
\numberwithin{equation}{section}
\def\simleq{\underset\sim<}
\def\simgeq{\underset\sim>}

\def\T{\text}
\def\1#1{\overline{#1}}
\def\2#1{\widetilde{#1}}
\def\3#1{\widehat{#1}}
\def\4#1{\mathbb{#1}}
\def\5#1{\frak{#1}}
\def\6#1{{\mathcal{#1}}}
\def\C{{\4C}}

\def\sumK{\underset{|K|=k-1}{{\sum}'}}
\def\sumJ{\underset{|J|=k}{{\sum}'}}

\def\Re{{\sf Re}\,}
\def\Im{{\sf Im}\,}

\def\phi{\varphi}

\def\Om{\Omega}

\newtheorem{Thm}{Theorem}[section]
\newtheorem{Cor}[Thm]{Corollary}
\newtheorem{Pro}[Thm]{Proposition}
\newtheorem{Lem}[Thm]{Lemma}
\theoremstyle{definition}\newtheorem{Def}[Thm]{Definition}
\theoremstyle{remark}
\newtheorem{Rem}[Thm]{Remark}
\newtheorem{Exa}[Thm]{Example}

\def\Label#1{\label{#1}}
\def\bl{\begin{Lem}}
\def\el{\end{Lem}}
\def\bp{\begin{Pro}}
\def\ep{\end{Pro}}
\def\bt{\begin{Thm}}
\def\et{\end{Thm}}
\def\bc{\begin{Cor}}
\def\ec{\end{Cor}}
\def\bd{\begin{Def}}
\def\ed{\end{Def}}
\def\br{\begin{Rem}}
\def\er{\end{Rem}}
\def\be{\begin{Exa}}
\def\ee{\end{Exa}}
\def\bpf{\begin{proof}}
\def\epf{\end{proof}}
\def\ben{\begin{enumerate}}
\def\een{\end{enumerate}}

\def\1alpha{[\frac1\alpha]}
\def\T{\text}

\def\C{{\Bbb C}}

\numberwithin{equation}{section}
\def\T{\text}

\newcommand{\bom}{\bar{\omega}}

\newcommand{\no}[1]{\|{#1}\|}
\newcommand{\NO}[1]{\|{#1}\|^2}

\newtheorem{theorem}{Theorem  }[section]
\newtheorem{definition}[theorem]{Definition }
\newtheorem{lemma}[theorem]{Lemma  }
\newtheorem{proposition}[theorem]{Proposition  }
\newtheorem{corollary}[theorem]{Corollary }
\newtheorem{example}[theorem]{\it Example }

\begin{document}
\title[The Diederich-Fornaess index and the global regularity...]{The Diederich-Fornaess index and the global regularity of the $\dib$-Neumann problem}
\author[ S.~Pinton and G.~Zampieri ]
{Stefano Pinton and Giuseppe Zampieri}
\address{Dipartimento di Matematica, Universit\`a di Padova, via 
Trieste 63, 35121 Padova, Italy}
\email{pinton@math.unipd.it,
zampieri@math.unipd.it}
\maketitle
\begin{abstract}
We describe along the guidelines of Kohn \cite{K99}, the constant $\mathcal E_s$ which is needed to  control  the commutator of a totally real vector field $T_{\mathcal E}$ with $\dib^*$ in order to have $H^s$ a-priori estimates for the Bergman projection $B_k$, $k\geq q-1$, on a smooth $q$-pseudoconvex domain $D\subset\subset \C^n$. This statement, not explicit in \cite{K99}, yields regularity results for $B_k$ in specific Sobolev degree $s$. Next, we refine the pseudodifferential calculus at the boundary in order to relate,  for a defining function $r$ of $D$, the operators $(T^+)^{-\frac\delta2}$ and $(-r)^{\frac\delta2}$. We are thus able to extend to general degree $k\ge0$ of $B_k$, the conclusion of \cite{K99} which only holds for $q=1$ and $k=0$: if for the Diederich-Fornaess index $\delta$ of $D$, we have $(1-\delta)^{\frac12}\le \mathcal E_s$, then $B_k$ is $H^s$-regular.

\end{abstract}
\section{Introduction}
\Label{s0}
The regularity of the Bergman projection $B_k$ over forms of degree $k\geq0$, as well as of the Neumann operator $N_k$ for $k\ge1$ on a pseudoconvex domain $D\subset\subset \C^n$,    has a long history. The first approach by Boas and Straube 
\cite{BS90}, \cite{BS91} consists in requiring, for any $\epsilon$, the existence of a totally real vector field $T_\epsilon
$, $|T_\epsilon|\sim1$, such that
\begin{equation}
\Label{0.0}
\Big|\langle \di r,[\di_{\bar z_i},T_\epsilon]\rangle\Big|\Big|_{bD}<\epsilon,\quad i=1,...,n,
\end{equation}
where $r$ is a defining function with $|\di r|=1$. 
This is referred to as ``good vector fields" condition. 
In other terms, we are requiring that all the coefficients of the $T_\epsilon$ components of $[\dib^*,T_\epsilon]$ are small
(modulo ``good" terms); cf. \cite{S10} Proposition 5.26. This can be weakened to a ``multiplier" condition for $[\dib^*,T_\epsilon]$. 
 Thus, the regularity of $B_k$, $k\ge0$  and $N_k$, $k\ge1$, is in fact related to the 
existence, for any $\epsilon$, of a totally real vector field $T_\epsilon$, with $|T_\epsilon|\sim1$, such that
\begin{equation}
\Label{0.1}
\NO{[\dib^*,T_\epsilon]u}<\epsilon Q_1(u,u)+c_\epsilon||u||,
\end{equation}
where $Q_1(u,u)=\NO{\dib u}_1+\NO{\dib^*u}_1$. 
Indeed, in \eqref{0.0} and \eqref{0.1} one can make the weaker assumtion that $T_\epsilon$ is ``approximately tangential", that is, $|T_\epsilon r|<\epsilon$; we refer for this point to the remarks after Theorem 5.22 of  \cite{S10}.
We deform the defining function $r$  to $r_\epsilon=g_\epsilon r$ and, accordingly, we deform the vector field $T=2\Im\frac{\sum_ir_{\bar i}\di_{z_i}}{\sum_i|r_i|^2}$ to $T_{g_\epsilon}=2\Im\frac{\sum_i(r_\epsilon)_{\bar i}\di_{z_i}}{\sum_i|(r_\epsilon)_i|^2}$.
The condition of approximate tangentiality turns into $|\Im g_\epsilon|<\epsilon$.
 These two deformations are related by $[\dib^*,T_{g_\epsilon}]\sim (\di\dib r_\epsilon\lefthalfcup\dib r_\epsilon) T_{g_\epsilon}$  modulo an error whose restriction to $bD$ belongs to $T^{1,0}bD\oplus T^{0,1}\C^n|_{bD}$; hence, the existence of $r_\epsilon$ such that
\begin{equation}
\Label{0.2}
|\di\dib r_\epsilon \lefthalfcup\dib r_\epsilon|\le \epsilon Q+c_\epsilon \Lambda^{-1},
\end{equation}
for $|\di r_\epsilon|\sim1$, implies \eqref{0.1}. (Here $\Lambda$ is the standard elliptic operator of order 1.) This is indeed the assumption under which   Straube proves  in \cite{S08} $H^s$-regularity for any $s$. In particular, this condition is fulfilled when there is a smooth defining function $r$ such that $\di\dib r|_{bD}\geq 0$; in this case one takes, for any $\epsilon$, $r_\epsilon=r$ in \eqref{0.2} and $T_\epsilon=T$   in \eqref{0.1} respectively (cf. the proof of Theorem~\ref{t1.2} below). Note that, historically, the conclusion was obtained, instead, through the ``good vector fields" condition. However how this follows from the fact that there exists $r$ which is plurisubharmonic on $bD$ is not immediate (Remark~\ref{r1.3} below). In any case, \eqref{0.0} calls into play a full family $\{T_\epsilon\}$ and the way of getting $T_\epsilon$ from the initial $T$ is involved.
In \cite{K99}, Kohn has given a quantitative result on regularity: he has specified, for given $s$,
and by allowing a full flexibility in the choice of $g$, not necessarily $g\sim1$, 
which is the constant $\mathcal E_{s,g}$ which is needed in \eqref{0.1} or \eqref{0.2} for $H^s$-regularity. This is not explicitly stated, but is entirely contained in \cite{K99} which, in turn, goes back to \cite{BS91}. If this is separated from the body of the paper, as we do in Theorem~\ref{t1.1}, and under an additional assumption of uniformity under exhaustion, it gives  $H^s$-estimates; this separation only requires minor modifications and yields a conclusion which naturally extends  to forms of any degree $k\ge q$ on $q$-pseudoconvex domains. 

It has been proved by Diederich-Fornaess in \cite{DF77} that every domain possesses an index $\delta$ with $0<\delta\le1$
such that $-(-r_\delta)^\delta$ is plurisubharmonic. Again, $r_\delta$ is in the form $r_\delta=g_\delta r$ for some $g_\delta$. On the other hand, it has been proved by Barret \cite{B92} that given a Sobolev index $s\searrow0$, one can find a domain $D$ in which $B_k$ fails $H^s$-regularity; according to \cite{DF77}, for these  domains, one has $\delta\searrow0$. So the relation between the index  of regularity $s$ and the Diederich-Fornaess index $\delta$ is an attractive problem. Indeed, what is explicitly stated by Kohn  and is by far the most interesting content of \cite{K99}, is the way of obtaining $\mathcal E_{s,g}$ out of $\delta$. This is described through  the estimate of the Levi form
\begin{equation*}
(-r_\delta)^{\frac\delta2}\Big|\di\dib(-(-r_\delta)^\delta)\lefthalfcup\dib r_\delta\Big|\simleq (1-\delta)^{\frac12}Q_{(-r_\delta)^{\frac{\delta}2}}.
\end{equation*}
(For an operator Op, such as Op$=(-r_\delta)^\delta$, we define $Q_{\T{Op}}$  by $Q_{\T{Op}}(u,u)=\NO{\T{Op}\dib u}+\NO{\T{Op}\dib^* u}$.)
In this estimate, one enjois the presence of the factor $(1-\delta)^{\frac12}$. When $(1-\delta)^{\frac12}\le \mathcal E_{s,g}$, one expects $s$-regularity by what has been said above, but this is not given for free because one encounters the unpleasant factor $(-r_\delta)^{\frac\delta2}$. It is well known that $(-r_\delta)^{\frac\delta2}\sim (T^+)^{-\frac\delta2}$ when the action is restricted to harmonic functions. For this reason, Kohn can prove regularity for the projection $B_0$ on $0$-forms, since this factorizes through the projection over harmonic functions. The main task of the present paper is to develop an accurate pseudodifferential calculus at the boundary which relates the action of $(-r_\delta)^{\frac\delta2}$ and  $ (T^+)^{-\frac\delta2}$ over general functions by describing the error terms by means of $\Delta$. In this way, when $(1-\delta)^{\frac12}\leq \mathcal E_{s,g}$, we get $H^s$-regularity of $B_k$ in general degree $k\geq0$ (resp. $k\ge q-1$) on a pseudoconvex (resp. $q$-pseudoconvex) domain.

Recent contribution to regularity of the Bergman projection by the method of the ``multiplier" is given by Straube in the already mentioned paper \cite{S08} and Herbig-McNeal \cite{HN06}; 
a combination of the ``multiplier" and ``potential" method  (inspired to the ``(P)-Property" by Catlin) is developed by Khanh \cite{Kh10} and Harrington \cite{H11}.
\vskip0.3cm
\noindent
{\it Acknowledgements.} We are grateful to Emil Straube for important advice.

\section{Weak $s$-compactness and $H^s$-regularity}
\Label{s1}
Let $D$ be a bounded smooth domain of $\C^n$ defined by $r<0$ for $\partial r\neq0$. We  modify the defining function as $gr$ for $g\in C^\infty$ and use the notation $r_g$ or $r^g$ for $gr$. We use the lower scripts $i$ and $\bar j$ to denote derivative in $\di_{z_i}$ and $\di_{\bar z_j}$ respectively and work with various vector fields such as 
\begin{equation}
\Label{1.0}
N_g=\frac{\sum_i r^g_{\bar i}\di_{z_i}}{\sum_i|r^g_i|^2},\qquad L^g_j=\di_{z_j}-r^g_jN_g,\qquad T_g=-i(N_g-\bar N_g).
\end{equation}
The $L^g_j$'s are complex-tangential; $T_g$ is the complementary real-tangential vector field.
 We consider an orthonormal basis $\bom_1,...,\bom_n$ of antiholomorphic 1-forms and general  forms $u$ of degree $k$, that is, expressions of type  $u=\sumJ u_J\bom_J$ where $J=j_1<...j_k$ are ordered multiindices and $\bom_J=\bom_1\wedge...\wedge \bom_k$. 
 We use the notations
 $$
\mathcal S=\T{Span}\{L_j^g,\,\di_{\bar z_j},\T{ for } j=1,...,n\}, \qquad
  Q_s(u,u)=\no{\dib u}^2_s+\no{\dib^* u}^2_s.
$$
We have (cf. \cite{BS91} p. 83) for $u\in C^\infty(\bar D)$,
\begin{equation}
\Label{1.1}
\no{Su}^2_{s-1}\simleq Q_{s-1}(u,u)+\no{u}_s\no{u}_{s-1} \quad\T{for any $S\in\mathcal S$}.
\end{equation}
Since $\mathcal S\oplus \C T_g=\C\otimes T\C^n$, then \eqref{1.1} implies
\begin{equation}
\Label{1.2}
\no{u}^2_s\simleq Q_{s-1}(u,u)+\no{T^s_gu}^2+\no{u}_s\no{u}_{s-1}.
\end{equation}
With the notation $\bar \theta_j:=-\frac1{\sum_i |r^g_i|^2}\sum_i r^g_{i\bar j}r^g_{\bar i}$, we define
\begin{equation}
\Label{1.3}
\begin{cases}
\bar{  \Theta}_gu=\sumK\sum_{ij}\Big(\bar\theta_j^gu_{iK}-\bar\theta^g_iu_{jK}\Big)+\T{error},
\\
\bar{  \Theta}_g^*u=\sumK\sum_{j}\theta_j^gu_{jK}+\T{error}.
\end{cases}
\end{equation}
We have the crucial commutation relation between $T_g$ and the Euclidean derivatives (\cite{K99} Lemma 3.33)
\begin{equation}
\Label{new}
[\di_{\bar z_j},T_g]=\bar \theta_jT_g\qquad \T{modulo $\mathcal S$}.
\end{equation}
This implies
\begin{equation}
\Label{1.4}
[\dib, T_g]=\bar{  \Theta}_gT_g\qquad\T{modulo $\mathcal S$}.
\end{equation}
As for the commutation of the adjoint $\dib^*$, we need a modification of $T_g$ which preserves the condition of membership to $D_{\dib^*}$. To this end, we define $\tilde T_g$ by
\begin{equation}
\Label{1.4,5}
(\tilde T_gu)_{jK}=T_gu_{jK}+\frac{r^g_{\bar j}}{\sum_i|r_{\bar i}|^2}\sum_i[T_g,r^g_i]u_{iK}.
\end{equation}
Thus $u\in D_{\dib^*}$ implies $\tilde T_g u\in D_{\bar \dib^*}$. 
Note that $\tilde T_g$ differs from $T_g$ by a $0$-order operator.
 With these preliminaries, \eqref{new} yields
\begin{equation}
\Label{1.5}
[\dib^*,\tilde T_g]=\bar{  \Theta}^*_g\tilde T_g\qquad\T{modulo $\mathcal S$}.
\end{equation}
\bd
\Label{d1.1}
Let $ s$ be a positive integer and let $1\leq q\leq n-1$. We say that $T^s_g$ well commutes with $\dib^*$ in degree $\geq q$ when
\begin{equation}
\Label{1.6}
\no{\bar \Theta^*_gu}^2\leq \mathcal E_{s,g}Q(u,u)+c_g\no{u}^2_{-1},\qquad\T{for any $u$ 
of degree $\geq q$},
\end{equation}
and for $\mathcal E_{s,g}\le c_1^2e^{-2c_2s\,\,\T{diam}^2\, D }\inf\Big(\frac1{|g|^s}\Big)^{-1}$ 
or, alternatively, for $\mathcal E_{s,g}\le c_1^2e^{-2c_2s\,\,\T{diam}^2D\sup(1+\frac{|g'|}{|g|})}$,
where $c_1$ is a small constant and $c_2$ is controlled by the $C^2$ norm of $r_g$.
\ed
We introduce the notion of $q$-pseudoconvexity of $D$; this consists in the requirement that, for the ordered eigenvalues $\lambda_1\le\lambda_2\le...\le\lambda_{n-1}$ of the Levi form $\di\dib r|_{\di r^\perp}$, we have $\underset{j=1}{\overset q\sum}\lambda_j\ge0$.
The basic estimates show that the complex Laplacian $\Box$ is invertible over $k$-forms for $k\ge q$. We denote by $N_k$ the inverse; we also denote by
$B_k:\,L^{2,k}\to L^{2,k}\cap \ker\dib$ the Bergman projection. Recall Kohn's formula $B_k=\T{Id}-\dib^*_{k+1}N_{k+1}\dib_k$. We say that $B_k$ is regular, resp. $s$-exactly regular, when it preserves $C^\infty$, respectively $H^s$, the $s$-Sobolev space.
\br
Assume that for any $s$ there is $r_g$ with $|\di r_g|\sim1$, that is $|g|\sim1$, such that $| \Theta^*_gu|\leq c_1e^{-c_2s\,\,\T{diam}^2\, D }$; then there is exact $s$-regularity for any $s$.
\er
We recall from \cite{BS90} that $s$-exact regularity of $N_k$ is equivalent to $s$-exact of the triplet $B_{k-1},\,\,B_k,\,\,B_{k+1}$. 
\bt
\Label{t1.1}
Let $D$ be $q$-peudoconvex and assume that for some $g$, $T_g^s$ well commutes with $\dib^*$ in degree $\geq q$.
Assume also that this property of good commutation holds, with a uniform constant $\mathcal E_{s,g}$, for a strongly $q$-pseudoconvex exhaustion of $D$.
 Then for any form $f\in H^s$  we have that $B_kf\in H^s$ and 
\begin{equation}
\Label{supernova}
\no{B_kf}_s\le c\no{f}_s,\quad\T{ for any $k\geq q-1$.}
\end{equation}
\et 
The proof is intimately related to \cite{BS91}. 
Formally, it follows the lines of \cite{K99} but also contains ideas taken from \cite{Kh10}.
\bpf
We first assume that we already know that $B_k$ is regular for any $k\geq q-1$ and prove \eqref{supernova} for a constant $c$ which only depends on \eqref{1.6}. In other terms, we show that \eqref{supernova} holds for $c$ if we knew from the beginning that it holds for some $c'>>c$.
We reason by induction.
An $n$ form is 0 at $b D$ ; thus $N_n$ ``gains two derivatives" by elliptic regularity of $\Box$ in the interior and hence $B_{n-1}$ is regular. We assume now that $B_k$ is $s$-regular and prove that the same is true for $B_{k-1}$. We use the notation $f$ for the test form in our proof; the notation $u$, which occurs in \eqref{1.6}, will be reserved to $\dib N_kf$. It suffices to estimate $\no{T^s_gB_{k-1}f}$ since, by \eqref{1.2}, this controls the full norm $\no{B_{k-1}f}_s$. We have
\begin{equation}
\Label{1.7}
\begin{split}
\no{T^s_gB_{k-1}f}^2&=\left(T^s_gB_{k-1}f,T^s_gf\right)-\left(T_g^sB_{k-1}f,T^s_g\dib^*N_{k}\dib f\right)
\\
&=\underset{(a)}{\underbrace{\left(T^s_gB_{k-1}f,T^s_gf\right)}}-\underset{(b)}{\underbrace{\left(T^{s*}_gT^s_g\dib B_{k-1}f,N_k\dib f\right)}}
\\
&{}\qquad-\underset{(c)}{\underbrace{\left([\dib,T^{s*}_gT^s_g]B_{k-1}f,N_k\dib f\right)}}.
\end{split}
\end{equation}
Now, $(a)\le sc\no{T^s_gB_{k-1}f}^2+lc\NO{T^s_gf}$, whereas $(b)=0$. 
The term which comes with small constant can be absorbed because we know a-priori that $\no{T^s_gB_{k-1}f}<\infty$. 
As for the last term, we replace $T_g^s$ by $\tilde T_g^s$ modulo an operator of order $s-1$, that we regard as an error term, describe the commutator in the left of (c) by $\bar{ \Theta}_g$ according to \eqref{1.4}, switch it to the right as  $ \bar\Theta_g^*$ and end up with
\begin{equation}
\Label{extra}
\begin{split}
|(c)|&\le \Big|\left(2s\bar{ \Theta}_g\tilde T^s_gB_{k-1}f,\tilde T_g^sN_k\dib f\right)\Big|+\T{error}
\\
&\le sc\NO{T^s_gB_{k-1}f}+lc\,s\NO{\bar \Theta^*_gT^s_gN_k\dib f}+\T{error}.
\end{split}
\end{equation}
The error includes terms in $(s-1)$-norm and terms in which derivatives belonging to $\mathcal S$ occur (cf. \eqref{1.1}). We use the hypothesis \eqref{1.5} under the choice $\mathcal E_{s,g}\le c_1^2c^{-2c_2s\,\,\T{diam}^2D}\sup\frac1{|g|^{2s}}$ and get,  with the notation $u=N_k\dib f$
\begin{equation}
\Label{1.8}
\begin{split}
\NO{ \bar\Theta^*_g\tilde T^s_gu}&\le\sup\frac1{|g|^{2s}}\NO{ \bar\Theta^*_g\tilde T^su}
\\&\le \mathcal E_{s,g}\sup\frac1{|g|^{2s}}Q(\tilde T^su,\tilde T^su)+\T{error}
\\&\le \mathcal E_{s,g}\sup\frac1{|g|^{2s}}\left(Q_{\tilde T^s}(u,u)+\NO{[\dib,\tilde T^s]u}+\NO{[\dib^*,\tilde T^s]u}\right)+\T{error}.
\end{split}
\end{equation}
(In case $\mathcal E_{s,g}\le c_1^2e^{-2c_2s\,\,\T{diam}^2D(1+\sup\frac{|g'|}{|g|})}$ we have not to replace $\tilde T^s_g$ by $\tilde T^s$ and, instead, use the estimate
\begin{equation}
\Label{g'/g}
|[\tilde T^s_g,\dib]v|\simleq c_2\sup(1+\frac{|g'|}{|g|})|\tilde T_gv|\quad\T{modulo $Sv$ for $S\in\mathcal S$,}
\end{equation}
and similarly for $\dib$ replaced by $\dib^*$; the proof will proceed similarly as below.)

Now,
$$
Q_{\tilde T^s}(u,u)\simleq \NO{T^sf}+\NO{T^sB_{k-1}f}+\T{error}.
$$
Next,
$$
\NO{[\dib,\tilde T^s]u}\leq c_2s^2\,\NO{T^sN_{k}\dib f}+\T{error}.
$$
We now observe that
\begin{equation}
\Label{1.-1}
\begin{split}
N_k\dib&=B_kN_k\dib(\T{Id}-B_{k-1})
\\
&=B_ke^{-\phi_s}N_{k,\phi_s}\dib e^{\phi_s}(\T{Id}-B_{k-1}),
\end{split}
\end{equation}
where $N_{k,\phi_s}$ is the $\dib$-Neumann operator weighted by $e^{-\phi_s}=e^{-c_2s|z|^2}$.
Since $[D^s,\dib]$ is an operator of degree $s$ with coefficients controlled by $sc_2$ for 
 $c_2\sim\no{r}_{C^2}$, then $N_{k,\phi_s}\dib$ is continuous in $H^s_{\phi_s}$
 with a continuity constant that we can assume to be unitary. 
  We use  that $c_2s^2\,e^{-2c_2s\,\,\T{diam}^2 D }\le \underset{z\in D}\inf e^{-2c_2s|z|^2}$ (for different $c_2$) in order to remove weights from the norms. We also use the inductive assumption that \eqref{supernova} holds for $B_k$. 
 In this way, we end up with 
\begin{equation}
\Label{nova}
\begin{split}
\mathcal E_{s,g}\sup\frac1{|g|^{2s}}c_2s^2\,\NO{T^s\dib N_kf}&\leq c_1^2\left(\NO{T^sf}+\NO{T^sB_{k-1}f}\right)+\T{error}
\\
&\le c_1^2(\NO{T^s_gf}+\NO{T^sB_{k-1}f})+\T{error},
\end{split}
\end{equation}
where the last inequality follows trivially from the fact that $T_g=\frac1g T$ for $\Big|\frac1g\Big|>>1$. 
Here, $\mathcal E_{s,g}$ takes care of $\sup\frac1{|g|^{2s}}$ and also of the constant which arises from removing weights owing to $\mathcal E_{s,g}\le c_1^2e^{-2sc_2\T{diam}^2 D }\sup \frac1{|g|^{2s}}$. Altogether, up to absorbable terms, $\NO{T^s_gB_{k-1}f}$ has been estimated by $lc\NO{T^s_gf}+\T{error}$. 
This  concludes the proof of Theorem~\ref{t1.1} if we are able to remove the assumption that we already know that \eqref{supernova} holds for some $c'>>c$. For this, we recall that we are assuming that there is a strongly $q$-pseudoconvex exhaustion $D_\rho\nearrow D$ which satisfies \eqref{1.6} uniformly with respect to $\rho$. We observe that \eqref{supernova} holds over $D_\rho$ for $c'=c'_\rho$. What has been proved above shows that it holds in fact with $c$ independent of $\rho$. 
Passing to the limit over $\rho$ we get \eqref{supernova} for $D$.

\epf
\bt
\Label{t1.2}
(Boas-Straube \cite{BS91})  
If there is a defining function $r$ such that for the eigenvalues $\mu_1\le...\le\mu_{n}$ of the full Levi form $\di\dib r$ (not restricted to $\di r^\perp$) we have $\underset{j=1}{\overset q\sum}\mu_j\ge0$, then, $B_k$ is exactly $H^s$-regular for any $s$ and any $k\geq q-1$.
\et
\bpf
The proof consists in proving that \eqref{1.6} holds for any $\epsilon$ and uniformly over an exhaustion of $D$. More precisely, we will show that for any $\epsilon$, for $\bar\Theta^*$ independent of $\epsilon$ (associated to a normalized defining function $r$), and for suitable $c_\epsilon$, we have
\begin{equation}
\Label{1.50}
\NO{  \bar\Theta^* u}\le\epsilon Q(u,u)+c_\epsilon|||{u}|||_{-1}\quad\T{for $u$ in degree $k\geq q$;}
\end{equation}
moreover, we will prove that \eqref{1.50} holds for a strongly $q$-pseudoconvex exhaustion. (Here, the triplet $|||\cdot|||$ denotes the tangential norm (cf. \cite{K79}).) 

\noindent
{\bf (a)} We begin by noticing that $\di\dib r+O(|r|)\T{Id}\ge0$ over $k$-forms for $k\ge q$.  We can then apply Cauchy-Schwartz inequality and get
\begin{equation}
\Label{1.10}
(r_{i\bar j})(u,\di r)\le (r_{i\bar j})(u,u)^{\frac12}+O(|r|^{\frac12})|u|.
\end{equation}
\\
\noindent
{\bf (b)} The Levi form is a ``$\frac12$-subelliptic multiplier" (cf. \cite{K79}), that is
\begin{equation}
\Label{1.11} |||\left((r_{i\bar j})(u,u)\right)^{\frac12}|||^2_{\frac12}\leq Q(u,u).
\end{equation}
This can be proved from the basic estimate 
$$
\int_D(r_{i\bar j})(Tu,u)dV\le Q(u,u),
$$
by using the microlocalization $T^+$ and its decomposition $T^+=(T^+)^{\frac12}(T^+)^{\frac12\,*}$. 
(Here $dV$ is the element of volume.)
Also,
by Sobolev interpolation, we have
\begin{equation}
\Label{1.12}
\begin{split}
\NO{(r_{i\bar j})(u,u)^{\frac12}}&\le \epsilon\NO{(r_{i\bar j})(u,u)^{\frac12}}_{\frac12}+c_\epsilon\NO{u}_{-1}
\\
&\le\epsilon Q(u,u)+c_\epsilon\NO{u}_{-1},
\end{split}
\end{equation}
where $c_\epsilon\sim \epsilon^{-1}||r||_{C^2}$. 
Finally, we estimate the norm of the last term in \eqref{1.10}. We have
\begin{equation}
\Label{1.13}
\begin{split}
\NO{(-r)^{\frac12}u}&\le\epsilon\NO{\zeta_\epsilon u}_0+\NO{(1-\zeta_\epsilon)u}_0
\\
&\le \epsilon\NO{u}_0+\NO{(1-\zeta_\epsilon)u}_0\simleq \epsilon Q(u,u)+\NO{(1-\zeta_\epsilon)u}_0,
\end{split}
\end{equation}
where $\zeta_\epsilon$ is a cut-off outside of the $\epsilon$-strip such that $|\dot \zeta_\epsilon|\simleq 
\frac1\epsilon$(with $\zeta_\epsilon\equiv1$ at $bD$). Moreover, we have
\begin{equation}
\Label{1.14}
\NO{(1-\zeta_\epsilon)u}_0\leq \epsilon^3\NO{(1-\zeta_\epsilon)u}_{1}+c_\epsilon\NO{(1-\zeta_\epsilon)u}_{-1},
\end{equation}
and,
\begin{equation}
\Label{1.15}
\begin{split}
\epsilon^3\NO{(1-\zeta_\epsilon)u}_1&\underset{(i)}\simleq  \epsilon^3 Q_0((1-\zeta_\epsilon)u,(1-\zeta_\epsilon)u)
\\
&\simleq\epsilon^3 Q_0(u,u)+\epsilon^3\NO{\dot\zeta_\epsilon u}_0
\\
&\simleq\epsilon^3 Q_0(u,u)+\epsilon^3\epsilon^{-2}\NO{u}_0
\\
&\underset{(ii)}\simleq 2\epsilon Q_0(u,u),
\end{split}
\end{equation}
where (i) is Garding inequality applied to $(1-\zeta_\epsilon)u|_{b D}\equiv0$ and (ii) follows from applying the basic estimate to $\NO{u}_0$. Putting together \eqref{1.10}--\eqref{1.15}, we get \eqref{1.50}.

\noindent
{\bf (c)} We consider the exhaustion of $D$ by the  domains $D_\rho$ defined by $r_\rho<0$ for $r_\rho=r+\rho e^{A|z|^2}$; by a suitable choice of $A$, these domains  are strongly $q$-pseudoconvex. 
We remark that $\di\dib r_\rho\simgeq -\no{r}_{C^2}|r_\rho|\,\T{Id}\ge -c|r_\rho|\,\T{Id}$ over $k$ forms for $k\geq q$. By Cauchy-Schwarz inequality we get
\begin{equation}
\Label{1.51}
(r_{i\bar j}^\rho)(u,\di r)\le (r^\rho_{i\bar j})(u,u)^{\frac12}+c|r_\rho|^{\frac12}|u|\quad\T{for $u$ of degree $k\geq q$}.
\end{equation}
The Levi form $(r^\rho_{i\bar j})$ is a $\frac12$-subelliptic multiplier (uniformly over $\rho$) and can be estimated as in (b) as well as the term with $O(|r_\rho|^{\frac12}) $.
Altogether, for fixed $\epsilon$ for any $\rho\le \rho_\epsilon$ and for $\bar\Theta^*_\rho$ associated to the definng function $r_\rho$, we have got 
\begin{equation}
\Label{1.55}
\NO{\bar\Theta^*_\rho u}\leq \epsilon Q_{D_\rho}(u,u)+c_\epsilon \NO{u}_{-1},
\end{equation}
uniformly with respect to $\rho$. Passing to the limit over $\rho$, yields \eqref{1.50}.

\epf

\bt
\Label{t1.3}
Let $D$ be $q$-pseudoconvex and assume that for any $\epsilon$ there is $|g_\epsilon|\sim 1$ such that
 \begin{equation}
  \Label{1.100}
 |  \bar\Theta^*_{g_\epsilon}(u)|\le\epsilon |u|^2\quad\T{on $bD$ for $u$ in degree $k\geq q$.}
\end{equation}
Then $B_k$ is exactly $H^s$-regular for any $s$ and any $k\geq q-1$.
\et

\bpf
\eqref{1.100} readily implies 
\begin{equation}
\Label{1.53}
\NO{  \bar\Theta^*_{g_\epsilon} u}\simleq\epsilon\NO{u}+\no{g_\epsilon r}_{C^2}\NO{(1-\zeta_\epsilon)u}\quad\T{for $u$ in degree $k\geq q$.}
\end{equation}
By plugging \eqref{1.100}   with the basic estimate $\NO{u}\simleq Q(u,u)$ 
and the Garding inequality $\no{g_\epsilon r}_{C^2}\NO{(1-\zeta_\epsilon)u}\simleq \epsilon Q(u,u)+c_\epsilon\NO{u}_{-1}$,
we get 
\begin{equation}
 \Label{1.8,5}
\NO{  \bar\Theta^*_{g_\epsilon} u}\simleq\epsilon Q(u,u)+c_\epsilon\NO{u}_{-1}\quad\T{for $u\in D_{\dib^*}$ of degree $k\geq q$.}
\end{equation}
This would give the $H^s$-regularity of $B_k$ if we were able to 
prove the stability of \eqref{1.100} under a strongly $q$-pseudoconvex exhaustion.
For this, we fix $\epsilon_o$ and $g_{\epsilon_o}r$ and approximate $D$ by $D_\rho$ defined by $g_{\epsilon_o}r+\rho e^{A|z|^2}$; for suitable fixed $A$, these are strongly $q$-pseudoconvex for any $\rho$. Also, if we rewrite $g_{\epsilon_o}r+\rho e^{A|z|^2}=g_{\epsilon_o,\rho}r_\rho$ for a normalized equation $r_\rho$ of $D_\rho$, we have
\begin{equation*}
\begin{cases}
g_{\epsilon_o,\rho}\underset{C^2}\to g_{\epsilon_o},
\\
r_\rho\underset{C^2}\to r.
\end{cases}
\end{equation*}
Hence
$$
\bar\Theta^*_{\epsilon_o,\rho}(u)\to \bar\Theta^*_{\epsilon_o}(u)\quad\T{uniformly over $u$}.
$$
 We then apply Theorem~\ref{t1.1} to each $\Om_\rho$ and by uniformity of the estimate with respect to $\rho$ we get that $B_kf$ belongs to $H^s$ and satisfies \eqref{supernova}. 

\epf
\br
\Label{r1.3}
We can give an alternative proof of Theorem~\ref{1.2} which uses Theorem~\ref{t1.3}. First, according to the lemma in \cite{BS91}, the existence of a plurisubharmonic defining function $r$ implies the vector fields condition \eqref{0.0}.
(If $r$ is only $q$-plurisubharmonic, \eqref{0.0} must be adpted by considering, similarly as in \eqref{1.100}, the action over  forms $u$ of degree $k\ge q$.)
 If we knew that the good vector fields $T_\epsilon$ are of type $T_{g_\epsilon}=-i(N_{g_\epsilon}-\bar N_{g_\epsilon})$, then, by \eqref{1.5} we would get \eqref{1.100} and reach the conclusion from Theorem~\ref{t1.3}. In the general case, by \cite{S10} Proposition 5.26, the condition of good vector fields implies \eqref{1.100}. (In that proposition, it is proved a generalization of \eqref{1.5}. For any tangential vector field $T_\epsilon$, not necessarily defined by \eqref{1.0}, if we denote by $g_\epsilon$ its $(N-\bar N)$-component, we have $[\dib^*,T_\epsilon]|_{bD}=\bar \Theta^*_{g_\epsilon}|_{bD}T_\epsilon$ modulo elliptic multipliers ($r$ and $\di r$) and $\frac12$-subelliptic multipliers ($\di\dib r$).)
\er
\br
\Label{r1.4}
We point out that in \cite{S08}, Straube proves that \eqref{1.8,5} suffices for exact $H^s$-regularity for any $s$. 
This requires heavy work since, differently from \eqref{1.100}, \eqref{1.8,5} is not tranferred from $\Om$ to $\Om_\rho$.
\er

\section{Pseudodifferential calculus at the boundary}
\Label{s2}
There is an important theory about the equivalence between  $(-r)^\sigma$ and microlocal powers $T^{-\sigma}$ over harmonic functions; we need to develop this theory and allow the action over general functions controlling errors coming from the Laplacian. In this  discussion, we do  not modify $r$ to $r_g$ and $T$ nor $T_g$.
Also, we still write $T$ but mean in fact its positive microlocalization $T^+$ which represents over $v^+$ the full elliptic standard operator $\Lambda$; for this reason, negative and fractionl powers of $T$ make sense.
We denote by $U$ a neighborhood of $bD$,
\bl
\Label{l2.1}
We have
\begin{equation}
\Label{2.1}
\no{(-r)^{\frac\delta2}r^\sigma T^\sigma v}\simleq lc\no{(-r)^{\frac\delta2}v}+sc\no{T^{-\frac\delta2}v}+sc\no{-rT^{-1-\frac\delta2}\Delta v}\quad\T{for any $v\in C^\infty(\bar D\cap U)$ and $\sigma>-\frac12$}.
\end{equation}
\el
This is a generalization of \cite{K99} Lemma~2.6 in which the extra terms with power $\frac\delta2$ do not occur.
\bpf
We have
\begin{equation*}
\begin{split}
\no{(-r)^{\frac\delta2}r^\sigma T^\sigma v}^2&=((-r)^{\delta+2\sigma}T^{2\sigma}v,v)
\\
&=-(\di_r(-r^{1+2\sigma+\delta})T^{2\sigma}v,v)
\\
&=2\Re((-r)^{1+2\sigma+\delta}\di_rT^{2\sigma}v,v)
\\
&\le lc\no{(-r)^{\frac\delta2}v}^2+sc\NO{(-r)^{1+2\sigma+\frac\delta2}\di_rT^{2\sigma+\frac\delta2-\frac\delta2}v}
\\
&\underset{(*)}\le lc\NO{(-r)^{\frac\delta2}v}+sc\NO{T^{-\frac\delta2}v}+sc\NO{-rT^{-1-\frac\delta2}\Delta v},
\end{split}
\end{equation*}
where $(*)$ follows from \cite{K99} (2.4) applied for $1+2\sigma+\frac\delta2>0$.

\epf
In \cite{K99} there is a result, Lemma~2.6, which applies to powers $>-\frac12$ of $-r$; we need a variant, still for negative powers, for terms involving $\di_rv$.
\bl
\Label{l2.2}
We have
\begin{equation}
\Label{2.2}
\no{(-r)^\sigma\di_rT^\sigma v}\simleq \no{v}+\no{rT^{-1}\Delta v}+\no{T^{-2}\Delta v},\quad v\in C^\infty(\bar D\cap U),\,\,\sigma>-\frac12.
\end{equation}
\el
\bpf
We have
$$
\Big(\di_r(-r)^{2\sigma+1}\di_rT^{2\sigma-2}v,\di_rv\Big)=-2\Re\Big((-r)^{2\sigma+1}\di_r^2T^{2\sigma-2}v,\di_rv\Big).
$$
Write $\di_r^2=\Delta+Tan\di_r+Tan^2\sim \Delta+T\di_r+T^2$.
For the three terms $\Delta$, $T^2$ and $T\di_r$, we have the three relations below, respectively
\begin{equation*}
\begin{cases}
\Big(T^{-2}\Delta v,(-r)^{2\sigma+1}T^{2\sigma}\di_rv\Big)\le \NO{T^{-2}\Delta v}+\NO{v},
\\
\begin{split}
\Big((-r)^{2\sigma+1}T^{2\sigma}v,\di_r v\Big)&=\Big((-r)^{2\sigma+1}T^{2\sigma+1}v,\di_rT^{-1}v\Big)
\\
&\simleq \NO{v}+\NO{-rT^{-1}\Delta v},
\end{split}
\\
\begin{split}
\Big((-r)^{2\sigma+1}\di_rT^{2\sigma-1}v,\di_rv\Big)&=\Big((-r)^{2\sigma+1}\di_rT^{(2\sigma+1)-1}v,\di_rT^{-1}v\Big)
\\
&\le \NO{v}+\NO{-rT^{-1}\Delta v},
\end{split}
\end{cases}
\end{equation*}
where the three inequalities come from Cauchy-Schwartz inequality combined with repeated use of \cite{K99} (2.4) (always under the choice $s=0$ with the notations therein). 
Finally, we have to estimate the error term 
\begin{equation}
\Label{2.3}
\Big((r)^{2\sigma+1}[\Delta,T^{2\sigma-2}]v,\di_rv\Big).
\end{equation}
We express the commutator in \eqref{2.3} as
$$
[\Delta,T^{2\sigma-2}]=T^{2\sigma-1}+\di_rT^{2\sigma-2}.
$$
Thus \eqref{2.3} splits into two terms to which the two inequalities below apply
\begin{equation*}
\begin{cases}
\begin{split}
\Big((-r)^{2\sigma+1}T^{2\sigma-1}v,\di_r v\Big)&=\Big((-r)^{2\sigma+1}T^{(2\sigma+1)-1}v,T^{-1}\di_rv\Big)
\\
&\le\NO{v}+\NO{-rT^{-1}\Delta v},
\end{split}
\\
\begin{split}
\Big((-r)^{2\sigma+1}\di_rT^{2\sigma-2}v,\di_rv\Big)&=\Big((-r)^{2\sigma+1}\di_rT^{2\sigma-1}v,T^{-1}\di_rv\Big)
\\
&\le \NO{v}+\NO{-rT^{-1}\Delta v}.
\end{split}
\end{cases}
\end{equation*}

\epf
We are ready for the main technical tool in interchanging powers of $-r$ and $T$.
\bp
\Label{p2.1}
We have
\begin{equation}
\Label{2.4}
\no{T^{-\frac\delta2}v}\simleq \no{(-r)^{\frac\delta2}v}+\no{-rT^{-1-\frac\delta2}\Delta v}+\no{(-r)^{\frac\delta2}T^{-2}\Delta v}.
\end{equation}
\ep
\bpf
We start from  \cite{K99} Lemma~2.11
\begin{equation*}
\begin{split}
\no{T^{-\frac\delta2}v}&\simleq\no{(-r_\delta)^{\frac\delta2}v}+\no{-rT^{-1-\frac\delta2}\Delta v }
\\
&+\sum_j\no{(-r_\delta)^{\frac\delta2}\di_{\bar z_j}T^{-1}v }.
\end{split}
\end{equation*}
Now, the first  and second terms in the right are good (in the right side of the estimate we wish to end with).
As for the last, we have
\begin{equation}
\Label{2.4,5}
\begin{split}
\sum_j\Big((-r_\delta)^{\frac\delta2}\di_{\bar z_j}T^{-1}v &,(-r_\delta)^{\frac\delta2}\di_{\bar z_j}T^{-1}v \Big)\leq\Big|\Big((-r_\delta)^{\frac\delta2}\Delta T^{-2}v ,(-r_\delta)^{\frac\delta2}v \Big)\Big|
\\
&+2\sum_j\Big|\Re\Big([(-r_\delta)^{\delta},\di_{ z_j}]\di_{\bar z_j}T^{-1}v,T^{-1} v\Big)\Big|.
\end{split}
\end{equation}
The first term in the right is estimated by
\begin{equation*}
\begin{split}
\Big|\Big((-r_\delta)^{\frac\delta2}\Delta T^{-2}v ,(-r_\delta)^{\frac\delta2}v \Big)\Big|&\le lc\no{(-r)^{\frac\delta2}v}+sc\no{(-r)^{\frac\delta2}(\di_r^2+\di_rT+T^2)T^{-2}v}
\\
&\le lc\no{(-r)^{\frac\delta2}v}+sc\Big(\no{(-r)^{\frac\delta2}T^{-2}\di_r^2v}+\no{(-r)^{\frac\delta2}\di_rT^{-1}v}\Big)
\\
&\le lc\no{(-r)^{\frac\delta2}v}+sc\Big(\no{(-r)^{\frac\delta2}T^{-2}\Delta v}+\no{T^{-\frac\delta2}v}\Big).
\end{split}
\end{equation*}
The second term in the right of \eqref{2.4,5} has the estimate
\begin{equation*}
\begin{split}
\Big|\Re\Big([(-r_\delta)^{\delta},\di_{ z_j}]T^{-1}v,T^{-1} v\Big)\Big|&\simleq\underset{(i)}{\underbrace{\Big|\Big((-r)^{-1+\delta+\epsilon}T^{-1+\frac\delta2+\epsilon}v,(-r)^{-\epsilon}T^{-\frac\delta2-\epsilon}v\Big)\Big|}}
\\&+\underset{(ii)}{\underbrace{\Big|\Big((-r)^{-1+\delta}\di_rT^{-1}v,T^{-1}v\Big)\Big|}}.
\end{split}
\end{equation*}
To estimate (i), we write $-1+\delta+\epsilon=\frac\delta2+(-1+\frac\delta2+\epsilon)=\frac\delta2+\sigma$ under the choice of $\epsilon>\frac12-\frac\delta2$ so that $-1+\frac\delta2+\epsilon>-\frac12$. We then apply Lemma~\ref{2.1} and get the estimate of (i)
$$
(i)\le lc\NO{(-r)^{\frac\delta2}v}+sc(\NO{T^{-\frac\delta2}v}+\no{-rT^{-1-\frac\delta2}\Delta v}).
$$
As for (ii) we have
\begin{equation*}
\begin{split}
(ii)&=\Big|\Big((-r)^{-1+\delta+(1-\frac\delta2-\epsilon)}\di_rT^{-1-\epsilon}v,(-r)^{-1+\frac\delta2+\epsilon}T^{-1+\epsilon}v\Big)\Big|
\\
&
\simleq sc\Big(\NO{T^{-\frac\delta2}v}+\no{-rT^{-1}\Delta v}+\no{T^{-2}\Delta v}\Big)
\\
&+lc\Big(\NO{(-r)^{\frac\delta2}v}+\no{-rT^{-1-\frac\delta2}\Delta v}\Big).
\end{split}
\end{equation*}
In fact, the  term with lc in the last line comes from Lemma~\ref{l2.1} applied for $\sigma=-1+\epsilon$ (which requires $\epsilon>\frac12$). The term with sc is estimated by the aid of Lemma~\ref{l2.2}
\begin{equation*}
\begin{split}
\no{(-r)^{-1+\delta+(1-\frac\delta2-\epsilon)}\di_rT^{-1-\epsilon}v}&=\no{(-r)^{\frac\delta2-\epsilon}\di_rT^{-1+(\frac\delta2-\epsilon)-\frac\delta2}v}
\\
&\underset{\eqref{2.2}}\simleq\no{T^{-\frac\delta2}v}+\no{-rT^{-1}\Delta v}+\no{T^{-2}\Delta v}.
\end{split}
\end{equation*}

\epf
We decompose now $v=v^{(h)}+v^{(0)}$ where $v^{(h)}$ is the harmonic extension and $v^{(0)}:=v-v^{(h)}$; note that $v^{(0)}|_{bD}\equiv0$.
We also recall the modification $\tilde T$ of $T$ defined by \eqref{1.4,5} and  designed to preserve $D_{\dib^*}$.
\bp
\Label{p2.2}
We have
\begin{equation}
\Label{2.5}
\no{[\tilde T^{s-\frac\delta2},\dib^*]v^{(h)}}\simleq \no{(-r)^{\frac\delta2}[\tilde T^s,\dib^*]v^{(h)}},\quad v\in C^\infty(\bar D\cap U).
\end{equation}
\ep
\br
In turn, by \eqref{1.5}, we have $[\tilde T^s,\dib^*]=s\bar\Theta \tilde T^s$, and therefore \eqref{2.5} implies
\begin{equation}
\Label{2.6}
\no{[\tilde T^{s-\frac\delta2},\dib^*]v^{(h)}}\simleq s\no{(-r)^{\frac\delta2}\bar\Theta \tilde T^sv^{(h)}}.
\end{equation}
\er
\bpf
In fact, Jacobi identity yields
$$
[\tilde T^s,\dib^*]=-\tilde T^{s-\frac\delta2}[\tilde T^{\frac\delta2},\dib^*]+\tilde T^{\frac\delta2}[\tilde T^{s-\frac\delta2},\dib^*]+[\tilde T^{s-\frac\delta2}[\tilde T^{\frac\delta2},\dib^*]].
$$
It follows
\begin{equation}
\Label{2.7}
\tilde T^{\frac\delta2} [\tilde T^{s-\frac\delta2} ,\dib^*]=[\tilde T^s ,\dib^*]+\tilde T^{s-\frac\delta2}[\tilde T^{\frac\delta2},\dib^*]-[\tilde T^{s-\frac\delta2}[\tilde T^{\frac\delta2},\dib^*]].
\end{equation}
We apply $\tilde T^{-\frac\delta2}$ to both sides of \eqref{2.7} and use Proposition~\ref{p2.1}. The conclusion will follow once we are able to show that 
$-r\tilde T^{-1-\frac\delta2}[\Delta,[\tilde T^s,\dib^*]]$ and $(-r)^{\frac\delta2}T^{2}[\Delta,T^s\dib^*]$ 
are error terms. In fact, we write
\begin{equation*}
\begin{split}
[\Delta,[\tilde T^s,\dib^*]]&=[\di_r^2+\di_rTan+Tan^2,Tan^s+\di_r Tan^{s-1}]
\\&=Tan^{s-1}+\di_rTan^s\simleq \tilde T^{s+1}+\di_r\tilde T^s\qquad\T{modulo $\mathcal S$}.
\end{split}
\end{equation*}
It follows
\begin{equation*}
\begin{cases}
\no{-rT^{-1\frac\delta2}[\Delta,[\tilde T^s,\dib^*]]v^{(h)}}\simleq \no{-rT^{s-\frac\delta2}v^{(h)}}+\no{-r\di_rT^{s-1-\frac\delta2}v^{(h)}}\underset{\T{ \cite{K99} (2.4) }}\simleq \no{T^{s-1-\frac\delta2}v^{(h)}},
\\
\no{(-r)^{\frac\delta2}T^{-2}[\Delta,[\tilde T^s,\dib^*]]v^{(h)}}\simleq \no{(-r)^{\frac\delta2}T^{s-1}v^{(h)}}+\no{(-r)^{\frac\delta2}\di_rT^{s-2}v^{(h)}}\underset{\T{ \cite{K99} (2.4)}}\simleq \no{T^{s-1-\frac\delta2}v^{(h)}}.
\end{cases}
\end{equation*}
\epf

\section{Non-smooth plurisubharmonic defing functions}
\Label{s3}
\bd
\Label{d3.1}
We say that $D$ has a Diederich-Fornaess index $\delta=\delta_s$ for $0<\delta\le1$  which controls the commutators of $\dib$ and $\dib^*$ with $D^s$ over forms in degree $k\ge q$, when there is $r_\delta=g_\delta r$ for $g_\delta\in C^\infty,\,\,g_\delta\neq0$, such that
\begin{equation}
\Label{3.1}
\begin{cases}
-(-r_\delta)^\delta\T{ is $q$-plurisubharmonic, that is, the sum of the first}
\\
\hskip1cm \T{$q$ eigenvalues of $\di\dib (-(-r_g)^\delta)$ is non-negative }
\\
(1-\delta_s)\le \mathcal E_{s,g},
\end{cases}
\end{equation}
where $\mathcal E_{s,g}$ can be chosen so that $\mathcal E_{s,g}\le c_1e^{-c_2s\,\,\T{diam}^2D}\sup\left(\frac1{|g|^s}\right)^{-1}$ or, alternatively, 
$\mathcal E_{s,g}\le c_1e^{-c_2s\,\,\T{diam}^2D\sup(1+\frac{|g'|}{|g|})}$.
\ed
Related to the above notion, is the condition
\begin{equation}
\Label{3.2}
\NO{(-r_\delta)^{\frac\delta2} \bar \Theta^*_gu}\le\mathcal E_{s,g}Q_{(-r_\delta)^{\frac\delta2}}(u,u),
\end{equation}
for $\delta\le 1$. 
\bt
\Label{t3.1}
If $D$  is $q$-pseudoconvex and has a Diederich-Fornaess index $\delta=\delta_s$ which controls the commutators of $(\dib,\dib^*)$ with $D^s$ in degree $k\ge q$, then $B_k$ is $s$-regular for $k\ge q$.
\et 
\br
The proof consists in showing that \eqref{3.1} implies \eqref{3.2} (points (a) and (b) below) and then showing that \eqref{3.2} implies the conclusion. Note that,
 when $\delta=1$, we have in fact the better conclusion contained in Theorem~\ref{t1.2}.
\er
\bpf
We decompose a form as $u=u^\tau+u^\nu$ where $u^\tau$ and $u^\nu$ are the tangential and normal component respectively. We have
\begin{equation}
\Label{3.3}
\begin{cases}
\NO{u^\nu}_1\le \sum_i\NO{\di_{\bar z_i}u^\nu}_0\simleq Q(u,u)
\\
\begin{split}
Q(u^\tau,u^\tau)&\le Q(u,u)+Q(u^\nu,u^\nu)
\\
&\simleq Q(u,u)+\NO{u^\nu}_1
\\
&\simleq Q(u,u).
\end{split}
\end{cases}
\end{equation}
Hence it suffices to prove \eqref{3.2}. The same conclusion also applies to the decompositin $u=u^{(h)}+u^{(0)}$ and, in general, to any decomposition in which either of the two terms is $0$ at $bD$.

\noindent
{\bf (a)} We have
\begin{equation}
\Label{3.4}
\Big|\di\dib r_\delta(u^\tau,\di r_\delta)\Big|\simleq (1-\delta)^{\frac12}(-r_\delta)^{-\frac\delta2}\Big(\di\dib(-(-r_\delta)^\delta)(u^\tau,u^\tau)\Big)^{\frac12}.
\end{equation}
To see it, we start from
$$
\di\dib(-(-r_\delta)^\delta)=\delta(-r_\delta)^{\delta-1}\di\dib r_\delta+(-r_\delta)^{\delta-2}\delta(1-\delta)\di r\otimes\dib r.
$$
In particular,
$$
\di\dib r_\delta=\frac1\delta(-r_\delta)^{1-\delta}\di\dib(-(-r_\delta)^\delta)-(-r_\delta)^{-1}(1-\delta)\di r\otimes\dib r.
$$
We suppose that $\delta$ is bounded away from $0$ and, indeed, that it approaches $1$; thus we disregard it in the following. We have
\begin{equation*}
\begin{split}
\di\dib & r_\delta(u,\di r_\delta)\sim(-r_\delta)^{1-\delta}\di\dib (-(-r_\delta)^\delta)(u,\di r_\delta)-(-r_\delta)^{-1}(1-\delta)\di r_\delta\otimes \dib r_\delta(u,\di r_\delta)
\\
&\le (-r_\delta)^{1-\delta}\Big(\di\dib(-(-r_\delta)^\delta)(u,u)\Big)^{\frac12}\Big((-r_\delta)^{-2+\delta}(1-\delta)|\di r_\delta|^2+O((-r_\delta)^{-1+\delta})\Big)^{\frac12}
\\
&\hskip8cm+(1-\delta)|\di r_\delta|^2(-r_\delta)^{-1}|\di r_\delta\cdot u|
\\
&\simleq\Big((1-\delta)^{\frac12}(-r_\delta)^{-\frac\delta2}+O(-r_\delta)^{\frac12-\frac\delta2})\Big)\Big(\di\dib(-(r_\delta)^\delta)(u,u)\Big)^{\frac12}+(1-\delta)|\di r_\delta|^2(-r_\delta)^{-1}|\di r_\delta\cdot u|.
\end{split}
\end{equation*}
Evaluation for $u=u^\tau$, yields \eqref{3.4}.

\noindent
{\bf (b)} We prove now \eqref{3.2} using the basic estimates. Generally, these apply to smooth plurisubharmonic defining functions. However, in \cite{K99}, Kohn has a version for H\"older continuous plurisubharmonic functions such as $-(-r_\delta)^\delta$. This implies the inequality $(*)$ below
\begin{equation}
\Label{3.5}
\begin{split}
\NO{(-r_\delta)^{\frac\delta2}\bar\Theta^*_gu^\tau}&\simeq \int_D(-r_\delta)^\delta\Big|\di\dib r_\delta(u^\tau,\di r_\delta)\Big|^2 dV
\\
&\underset{\eqref{3.4}}\simleq (1-\delta)\int_D\di\dib(-(-r_\delta)^\delta)(u^\tau,u^\tau) dV
\\
&\underset{(*)}\simleq (1-\delta)Q_{(-r_\delta)^{\frac\delta2}}(u^\tau,u^\tau)
\\
&\underset{\eqref{3.1}}\simleq \mathcal E_{s,g}Q_{(-r_\delta)^{\frac\delta2}}(u^\tau,u^\tau).
\end{split}
\end{equation}
This proves \eqref{3.2}

\noindent
{\bf (c)}  We are therefore in the same situation as in Definition~\ref{d1.1} apart from the term $(-r_\delta)^\delta$ which occurs in the integral in the left of \eqref{3.5} and in $Q_{(-r_\delta)^{\frac\delta2}}$. 
As above, we continue to write $T$ but take in fact its positive microlocalization $T^+$ which represents the full action of $\Lambda$ over $u^+$.
To carry on the proof, we suppose from now on that $f\in C^\infty(\bar D)$ and that $B_k$ is $H^s$ regular for some continuity constant $c'$; we prove that this implies continuity for a constant $c$ which is solely related to the constants which occur in \cite{3.1}. An exhaustion by  domains endowed with $H^s$-regular projections $B_k,\,\,k\ge q$, will be discussed only at the end. 
We start from \eqref{1.7}
\begin{equation}
\Label{3.6}
\begin{split}
\no{T_g^{s-\frac\delta2}B_{k-1}f}&\simleq sc\NO{T^{s-\frac\delta2}_gB_{k-1}f}+lc\no{T^{s-\frac\delta2}_gf}
\\
&+lc\no{[\dib^*,T^{s-\frac\delta2}_g]N_k\dib f}.
\end{split}
\end{equation}
At this point, we need to convert $T^{s-\frac\delta2}_g$  into $(-r_\delta)^{\frac\delta2} T^s_g$ in the last term of \eqref{3.6} in order to enjoy \eqref{3.2}. We also replace $N_k\dib f$ by $(N_k\dib f)^{(h)}$ where the supscript $(h)$ denotes the harmonic extension. 
We apply the crucial estimate \eqref{2.5} to the last term in \eqref{3.6}, regard as errors the terms which come in $(s-1)$-norm or in which vector fields of $\mathcal S$ occur, and get
\begin{equation}
\Label{3.7}
\begin{split}
||[\dib^*,\tilde T^{s-\frac\delta2}_g]&(\dib N_kf)^{(h)}\underset{\eqref{2.5}}\leq\NO{(-r_\delta)^{\frac\delta2}[\tilde T^s_g,\dib^*](\dib N_kf)^{(h)}}
\\
&\simleq s^2\NO{(-r_\delta)^{\frac\delta2}\bar\Theta^*_g\tilde T^s_g(\dib N_kf)^{(h)}}+\T{error}
\\
&\simleq s^2\NO{(-r_\delta)^{\frac\delta2}\bar\Theta^*_g\tilde T^s_g(\dib N_kf)^{(h)\,\tau}}+\T{error}
\\
&\simleq s^2\sup\frac1{|g|^{2s}}\NO{(-r_\delta)^{\frac\delta2}\bar\Theta^*_g\tilde T^s(\dib N_kf)^{\tau}}+\mathcal E^{(0)}+\T{error}
\\
&\simleq s^2\mathcal E_{s,g}\sup\frac1{|g|^{2s}}\Big(Q_{(-r_\delta)^{\frac\delta2}\tilde T^s}((\dib N_kf)^\tau,(\dib N_kf)^\tau)
\\
&\hskip1cm +\NO{(-r_\delta)^{\frac\delta2}[\dib,\tilde T^s](\dib N_kf)^\tau}+\NO{(-r_\delta)^{\frac\delta2}[\dib^*,\tilde T^s](\dib N_kf)^\tau}\Big)+\mathcal E^{(0)}+\T{error}
\\
&\simleq s^2\mathcal E_{s,g}\sup\frac1{|g|^{2s}}\Big(Q_{(-r_\delta)^{\frac\delta2}\tilde T^s}(\dib N_kf,\dib N_kf)
\\
&\hskip1cm +\NO{(-r_\delta)^{\frac\delta2}[\dib,\tilde T^s]\dib N_kf}+\NO{(-r_\delta)^{\frac\delta2}[\dib^*,\tilde T^s]\dib N_kf}\Big)+\mathcal E^{(0)}+\T{error}
\\
&\simleq s^2\mathcal E_{s,g}\sup\frac 1{|g|^{2s}}\Big(\NO{(-r_\delta)^{\frac\delta2} T^s\dib^*\dib N_kf}+c_2s^2\NO{(-r_\delta)^{\frac\delta2} T^s\dib N_kf}\Big)+\mathcal E^{(0)}+\T{error}
\\
&\simleq s^2\mathcal E_{s,g}\sup \frac1{|g|^{2s}}\Big(\NO{(-r_\delta)^{\frac\delta2} T^s\dib^*\dib N_kf}+e^{2c_2\,\,s\,\,\T{diam}^2D}c_2s^2\NO{(-r_\delta)^{\frac\delta2} T^s\dib^*\dib N_kf}\Big)+\mathcal E^{(0)}+\T{error}
\\
&\underset{\eqref{3.1}}\simleq sc\NO{(-r_\delta)^{\frac\delta2}T^s\dib^*\dib N_kf}+\mathcal E^{(0)}+\T{error},
\end{split}
\end{equation}
where we have used the notation $\mathcal E^{(0)}:=\NO{(-r_\delta)^{\frac\delta2}\bar\Theta^*_g\tilde T^s_g(\dib N_kf)^{(0)\,\tau}}$.
Here in \eqref{3.1} we have chosen the first alternative $s^2\mathcal E_{s,g}e^{c_2s\,\,\T{diam}^2D}\sup\left(\frac1{|g|^s}\right)\le c_1=sc$ (for a new $c_2$). (The other alternative
$\mathcal E_{s,g}e^{c_2s\,\,\T{diam}^2D\sup(1+\frac{|g'|}{|g|}}\le c_1=sc$ can be handled similarly as in Theorem~\ref{t1.1}  without replacing $T_g$ by $T$. 
It is at this point, where the continuity of $B_k$ in $H^s$, not just in $C^\infty$, is needed; in fact, in formula \eqref{1.-1} $N_{\phi_s}$ is  $H^s$, not $C^\infty$, continuous.
We have to reconvert now $(-r_\delta)^{\frac\delta2}$ into $T^{-\frac\delta2}$. We first suppose that we had started from $f^{(h)}$ and wished to prove the regularity for $B_{k-1}f^{(h)}$. We have
\begin{equation*}
\no{(-r_\delta)^{\frac\delta2}T^s\dib^*\dib N_k(f^{(h)})}\underset{\T{  \cite{K99} (2.4)}}\simleq\underset{(i)}{\underbrace{ \no{T^{s-\frac\delta2}\dib^*\dib N_kf^{(h)}}}}+\underset{(ii)}{\underbrace{\no{-rT^{s-\frac\delta2-1}\Delta\dib^*\dib N_kf^{(h)}}}}.
\end{equation*}
Now, 
$$
(i)\simleq \NO{T^{s-\frac\delta2}f^{(h)}}+\NO{T^{s-\frac\delta2}B_{k-1}f^{(h)}},
$$
where the first term in the right is good and the second can be absorbed since it comes, inside \eqref{3.7}, with sc. As for (ii),
\begin{equation*}
\begin{split}
(ii)&=\no{-rT^{s-\frac\delta2-1}(\dib^*\dib+\dib\dib^*)\dib^*\dib N_kf^{(h)}}+\T{error}
\\
&=\no{-rT^{s-\frac\delta2-1}(\dib^*\dib(\dib\dib^*+\dib^*\dib )N_kf^{(h)}}+\T{error}
\\
&=\no{-rT^{s\frac\delta2-1}\dib^*\dib f^{(h)}}+\T{error}.
\end{split}
\end{equation*}
We have
\begin{equation*}
\begin{cases}
\dib^*\dib=Tan^2+\di_rTan+\di_r^2\sim T^2+\di_rT+\di_r^2,
\\
\di_r^2=\Delta+Tan^2+\di_rTan\sim \Delta+T^2+\di_rT,
\end{cases}
\end{equation*}
which implies
$$
\dib^*\dib\sim T^2+\di_rT+\Delta.
$$
It follows
\begin{equation}
\Label{3.8}
\begin{split}
\no{-rT^{s-\frac\delta2-1}\dib^*\dib f^{(h)}}&=\no{-rT^{s-\frac\delta2-1}(T^2+\di_rT+\Delta)f^{(h)}}
\\
&\leq \no{-rT^{s-\frac\delta2+1}f^{(h)}}+\no{-rT^{s-\frac\delta2}\di_rf^{(h)}}
\\
&\underset{\T{\cite{K99} (2.4)}}\simleq \no{T^{s-\frac\delta2}f^{(h)}},
\end{split}
\end{equation}
which is good. As for the term $f^{(0)}$, the regularity of $B_{k-1}f^{(0)}$ follows readily, without using the machinery (a)--(c) above, from elliptic regularity
\begin{equation}
\Label{3.9}
\no{T^sN_{k-1}f^{(0)}}\simleq \no{T^{s-2}f^{(0)}}.
\end{equation}
(Note that $N_{k-1}$ makes sense even for $k-1=0$ when acting on $f^{(0)}|_{bD}\equiv0$ because $\Box$ is, under this restriction, invertible.)

We pass to the term which has been omitted in the estimate of $\bar\Theta^*_g$, that is, $\mathcal E^{(0)}$. The use of elliptic regularity is different here and applies to $(\dib N_kf)^{(0)}$ instead of $f^{(0)}$; it then passes though $Q$ instead of $\Box$ and through Boas-Straube formula. We have
\begin{equation}
\Label{3.10}
\begin{split}
||(-r_\delta)^{\frac\delta2}&\bar\Theta^*_g\tilde T^s_g(\dib N_kf)^{(0)\,\tau}||^2
\simleq \sup\frac1{|g|^{2s}}\NO{(-r_\delta)^{\frac\delta2}\bar\Theta^*_g\tilde T^s(\dib N_kf)^{(0)\tau}}
\\
&\simleq \mathcal E_{s,g}\sup\frac1{|g|^{2s}}\Big(Q_{(-r_\delta)^{\frac\delta2}\tilde T^s}((\dib N_kf)^{(0)\tau},(\dib N_kf)^{(0)\tau)})+\T{error}
\\
&\hskip1cm +\NO{(-r_\delta)^{\frac\delta2}[\dib,\tilde T^s](\dib N_kf)^{(0)\tau}}+\NO{(-r_\delta)^{\frac\delta2}[\dib^*,\tilde T^s](\dib N_kf)^{(0)\tau}}\Big)+\T{error}
\\
&\simleq \mathcal E_{s,g}\sup\frac1{|g|^{2s}}\Big(Q_{(-r_\delta)^{\frac\delta2}\tilde T^s}(\dib N_kf,\dib N_kf)+\T{error}
\\
&\hskip1cm +\NO{(-r_\delta)^{\frac\delta2}[\dib,\tilde T^s]\dib N_kf}+\NO{(-r_\delta)^{\frac\delta2}[\dib^*,\tilde T^s]\dib N_kf}\Big)+\T{error}
\end{split}
\end{equation}
This is the same as \eqref{3.7} with the advantage that in the last line the Sobolev indices have decreased by $-1$ since terms with superscript $(0)$ vanish at $bD$; these are therefore error terms. Also there remain to control $\no{T^{-\frac\delta2}\bar\Theta^*_g\tilde T^s_g(\dib N_kf)^{(0)}}$ and $\no{-r^{\frac\delta2}\bar\Theta^*_g\tilde T^s_g(\dib N_kf)^\nu}$; but these are controlled by elliptic regularity as in \eqref{3.10}.
Summarizing up, we have proved that 
for a suitable $c$, only related to the constants in \eqref{3.1}, we have 
\begin{equation}
\Label{final}
\no{B_k f}_s\le c\no{f}_s
\end{equation}
 if we knew that it holds for some $c'>>c$.
We show now that we can exhaust $D$ by  domains $D_\rho$  endowed with continuous projections $B_k,\,\, k\geq q-1$ for some $c'$ and which inherit the assumption of Theorem~\ref{t3.1} with uniform constants with respect to $\rho$. For this, we define $D_\rho=\{z:\,r_\delta(z)+\rho<0\}$. We first notice that, $bD_\rho$ being also defined by $-(-r_\delta)^\delta+\rho^\delta<0$, it has a smooth $q$-plurisubharmonic defining function.
Hence, by Theorem~\ref{t1.2}, $B_k$ is $H^s$-regular for any $k\ge q-1$. Coming back to the initial defining function $r_\delta+\rho$, this satisfies $\di\dib (-(-r_\delta-\rho)^\delta)\ge \di\dib(-(-r_\delta)^\delta$; thus the Diederich-Fornaess index of $D_\rho$ is $\ge \delta$. Also, if for the new boundary we rewrite $r_\delta+\rho=g_{\delta,\rho}r_\delta$, for a normalized equation $r_\rho$ of $D_\rho$, and if $\mathcal E_{s,g,\rho}$ are the constants which occur in \eqref{3.1}, then
\begin{equation*}
\begin{cases}
g_{\delta,\rho}\underset{C^2}\to g_\delta,
\\
\mathcal E_{s,g,\rho}\underset{C^2}\to \mathcal E_{s,g}.
\end{cases}
\end{equation*}
Thus, the estimate \eqref{final} passes from the $D_\rho$'s (in which it has been proved thanks to the regularity of the $B_k$ (for a different $c'$)) to the initial domain $D$.

The proof is complete.

\epf

\end{document}